\documentclass[12pt]{amsart}
\usepackage{amssymb}
\usepackage{amsmath}
\usepackage[dvips]{graphicx}
\usepackage{psfrag}

\newcommand{\be}{\begin{equation}}
\newcommand{\ee}{\end{equation}}
\newcommand{\fd}{Ferrers diagram }
\newcommand{\fds}{Ferrers diagrams }
\newcommand{\fp}{Ferrers wicket }

\newtheorem{theo}{Theorem}

\newtheorem{defn}{Definition}

\title[Gated and Wicketed Ferrers diagrams ]{The Half-Perimeter Generating Function of Gated and Wicketed Ferrers diagrams}
\author{Arvind Ayyer}
\address{Arvind Ayyer\\
Department of Physics\\
136 Frelinghuysen Rd\\
Piscataway, NJ 08854.}
\email{ayyer@physics.rutgers.edu}
\date{\today}

\begin{document}

\begin{abstract}
We show that the half-perimeter generating functions for the number of Wicketed and Gated Ferrers diagrams is algebraic. Furthermore, the generating function of the Wicketed Ferrers diagrams is closely related to the generating function of the Catalan numbers. The methodology of the experimentation as well as the proof is the umbral transfer matrix method.
\end{abstract}

\maketitle

\section{Introduction}
Motivated by the recent interest in the enumeration of staircase polygons with a single staircase puncture and conjectures of a holonomic solution \cite{gutt}, we investigate the simpler problem of Ferrers diagrams with Ferrers punctures of different kinds - wickets and gates. Since Ferrers diagrams form a subset of staircase polygons, we simple-mindedly expect a holonomic solution here too. Fortunately for us, this na\"{\i}ve expectation turns out to be true and the generating function in both cases is not only holonomic, but also algebraic and moreover, the degree of the algebraic equation satisfied by the generating function is two! 

The umbral transfer matrix method \cite{zeil} is a technique to calculate terms in the series expansion of generating functions. For a given combinatorial building, the umbral operator is essentially the architectural plan for the structure. At any stage of the construction, it tells the builder exactly how to proceed from there on. The power of the method is twofold: it helps in generating terms in the sequence and secondly, once an ansatz is in place, it is very easy to prove or disprove the ansatz. It is in that spirit that the proofs here must be read, namely as simple exercises in algebra. The techniques here are more important than the proofs. 

We also emphasize the experimental nature of the paper. The proofs here are completely computerizable. For this problem, the solution is simple enough for everything to be done by hand. For more complex problems of this type, however, pen-and-paper calculations would be far too long and error-prone to be efficient. This method can be generalized to more complex problems such as the conjecture of punctured staircase polygons \cite{gutt}.

The plan of the paper is as follows: We first introduce the method by applying it to the simple case of standard Ferrers diagrams. We then go on to apply it to the case of interest deriving the various umbral operators that arise and proving the main theorems. Finally, we offer hope for a bijective proof of the main theorem.

\section{Standard Ferrers Diagrams}
\begin{defn}
A \emph{Ferrers Diagram} is a collection of $n$ rows of blocks, the $i$th row of which contains $m_i$ blocks, the first row at the bottom and the last row on top. All the rows are left aligned and such that if $1 \leq i<j\leq n$ then $m_1 \leq m_i \leq m_j \leq m_n$.
\end{defn}

Since a picture is worth an arbitrary number of words\footnote{More precisely, given a positive integer $N$, there exists a concept, a sentence consisting of $N$ words and a picture such that the picture describes the concept better than the sentence.}, and because we want our convention to be clear, we augment the definition with an example.

\begin{figure}[h!]
\includegraphics[width=6cm]{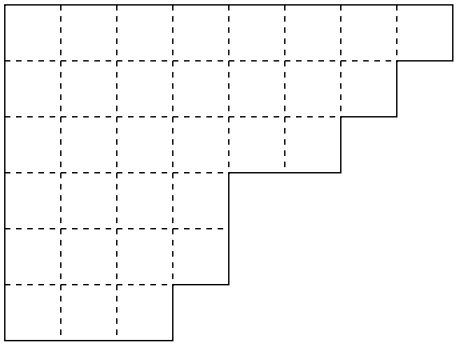}
\caption{A \fd with six rows, $m_1 = 3, m_2 = m_3 = 4, m_4 =6, m_5 = 7, m_6 = 8$. } \label{fig:fdeg}
\end{figure}

We first illustrate the methodology by applying it to the simple and almost trivial case of the half-perimeter generating function of usual Ferrers diagrams. The basic idea is to assign a catalytic variable to keep track of the width of the topmost segment and construct the evolution operator as a function of that variable.

\begin{figure}[h!]
\psfrag{a}{$a$}
\psfrag{b}{$b$}
\includegraphics[width=6cm]{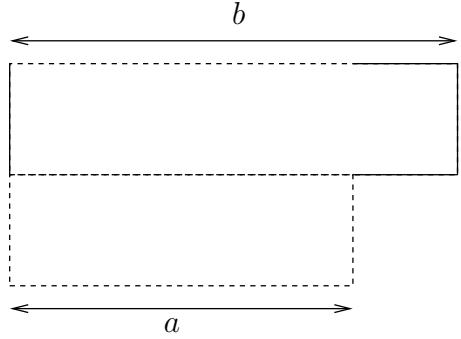}
\caption{The umbral evolution operator for Ferrers diagrams. The darkened lines contribute to the perimeter.} \label{fig:fdumb}
\end{figure}

Suppose, the top of the \fd has height $a$ as shown. Then, we construct all possible legal \fds which can be built from it. Let $x$ be the catalytic variable counting the width of the topmost segment and $t$ be the half-perimeter. The umbral evolution acts by adding the contributions of the top width and the extra half-perimeter for all such legal \fds,
\be
\begin{split}
x^a  \mapsto & \sum_{b=a}^\infty x^b t^{b-a+1} \\
& = \frac{x^a t}{1-x t}.
\end{split}
\ee

This gives the evolution for the monomial $x^a$. Since we need the operator evolution, we need the evolution for a formal power series $p(x)$:
\be \label{outpunc}
U_0(p(x)) = \frac{t p(x)}{1-xt}.
\ee

Now we need initial conditions. Obviously, they correspond to \fds of height one. Thus, as a formal power series, the initial condition is
\be
\begin{split}
I(x) & = \sum_{a=1}^\infty x^a t^{a+1} \\
& = \frac{x t^2}{1-xt}.
\end{split}
\ee

Then the generating function of \fds $F(x)$ satisfies the equation
\be
\begin{split}
F(x) & = I(x) + U_0(F(x)) \\
& = \frac{x t^2}{1-xt} + \frac{t F(x)}{1-xt},
\end{split}
\ee

which has the simple rational solution
\be
F(x) = \frac{x t^2}{1-t-xt},
\ee

Setting $x=1$ gives the well-known formula
\be \label{fdgf}
F(x) = \frac{t^2}{1-2t},
\ee
the expansion of which gives the terms in A000079 \cite{sloane}.

\section{Gated and Wicketed \fds}
\begin{defn}
A \emph{Gated \fd} is a \fd from which  another \fd with the same conventions is removed from the top.
\end{defn}

\begin{defn}
A \emph{Wicketed \fd} is a \fd from which another \fd with the same conventions is removed from the interior.
\end{defn}

\begin{figure}[h!]
\includegraphics[width=6cm]{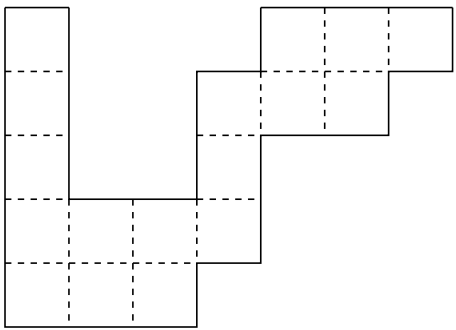}
\caption{A Gated \fd with five rows containing a three-rowed Ferrers gate.} \label{fig:fptopeg}
\end{figure}

\begin{figure}[h!]
\includegraphics[width=6cm]{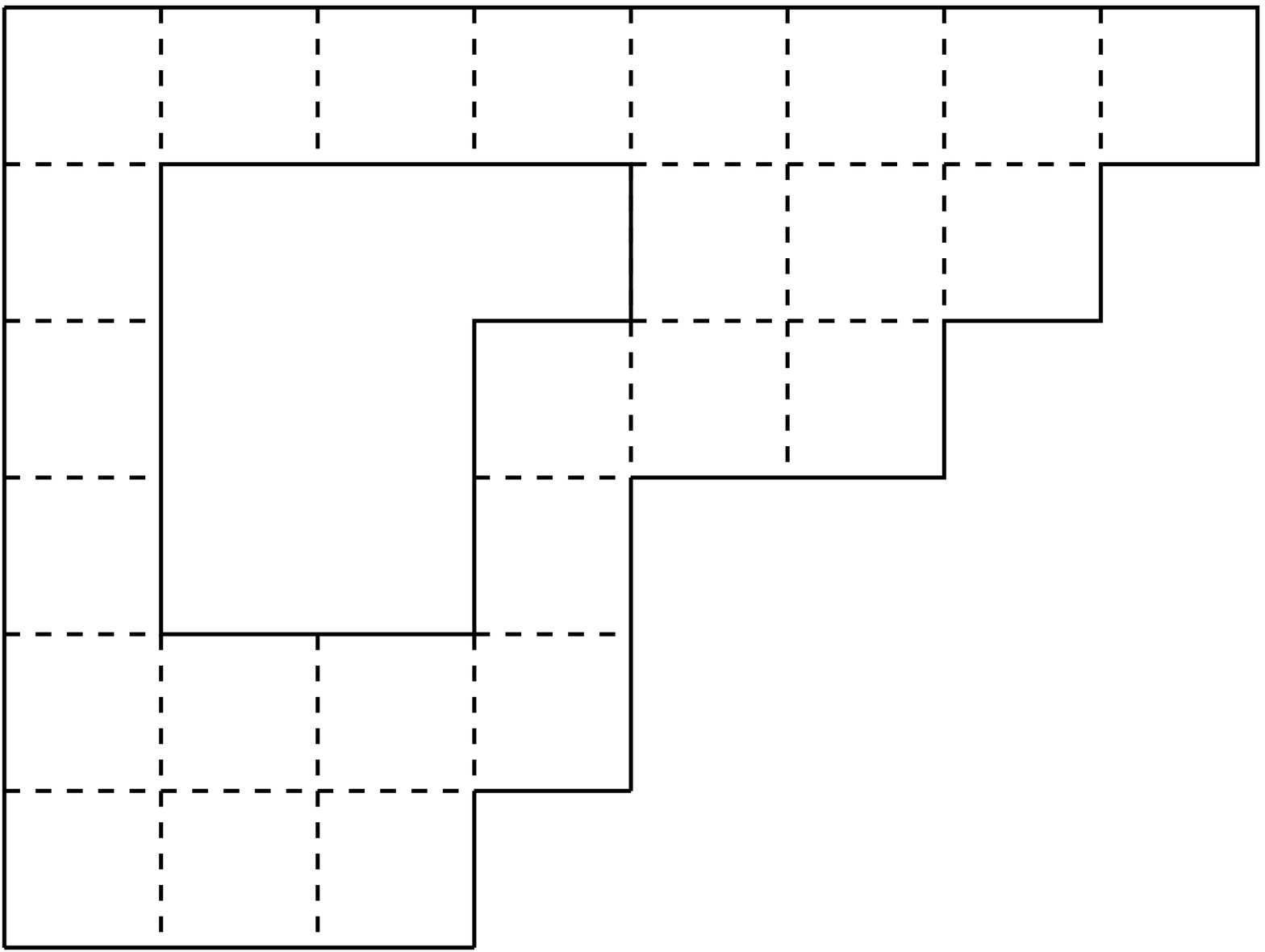}
\caption{A Wicketed \fd with six rows containing a three-rowed \fp.} \label{fig:fdfpeg}
\end{figure}

The terminology is inspired by the figures. A Gated Ferrers diagram like the one in Figure \ref{fig:fptopeg}, when inverted, looks as if a door or a gate is cut out from the \fd because the cut is at the bottom. A Wicketed \fd like the one in Figure \ref{fig:fdfpeg} looks as if there is a wicket or window cut out because the cut is in the interior. 

We use the umbral transfer matrix method again to count these objects. Unlike the simple case before, however, we need three and five different umbral operators for gated and wicketed diagrams respectively. Further complications arise from the fact that there are different number and types of catalytic variables in these different regions. 

\subsection{Construction of the gate/wicket}
Suppose that the construction before the wicket is complete. In other words, we have a Ferrers diagram so far. We now begin creating the wicket. We take an object with one catalytic variable $x$ counting the top width and add to it a different object with three catalytic variables $x_1,x_2,x_3$ marking the start of the wicket, the width of the wicket and the remaining width respectively also at the top.

\begin{figure}[h!]
\psfrag{a}{$a$}
\psfrag{b1}{$b_1$}
\psfrag{b2}{$b_2$}
\psfrag{b3}{$b_3$}
\includegraphics[width=6cm]{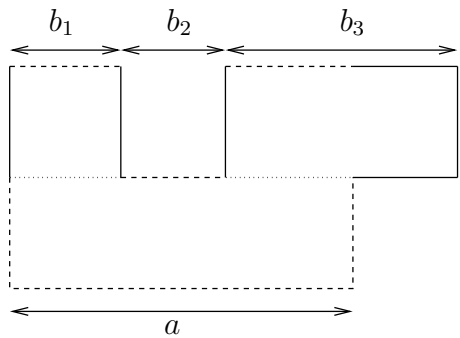}
\caption{The umbral evolution operator for creating the gate or wicket. The darkened lines contribute to the perimeter.} \label{fig:startpunc}
\end{figure}

Since we want the wicket to be strictly within the diagram, we can see that $b_1<a-1$ and $b_1+b_2<a$. The extra contribution to the half-perimeter can be seen by looking at the darkened lines in Figure \ref{fig:startpunc}. Thus the evolution operator $U_1$ acts on the monomial by
\be
\begin{split}
x^a \mapsto  & \sum_{b_1=1}^{a-2} \sum_{b_2=1}^{a-1-b_1} \sum_{b_3=a-b_1-b_2}^\infty x_1^{b_1} x_2^{b_2} x_3^{b_3} t^{2+b_1+b_2+b_3-a} \\
& = \frac{x_1^{a-1} x_2^2 x_3 t^2}{(x_1-x_2) (1-x_3 t) (x_2-x_3)}-\frac{x_1^{a-1} x_2 x_3^2 t^2 }{(x_1-x_3)(1-x_3 t)(x_2-x_3)}\\
&-\frac{x_1 x_2^a x_3 t^2}{(x_1-x_2) (1-x_3 t)(x_2-x_3)}+\frac{x_1 x_2 x_3^a t^2 }{(x_1-x_3)(1-x_3 t)(x_2-x_3)}
\end{split}
\ee

Therefore, $U_1$ acts on formal power series as follows:
\be \label{startpunc}
\begin{split}
&U_1(p(x)) = \frac{x_2 x_3 t^2 p(x_1)}{(x_1-x_2) (x_1-x_3) (1-x_3 t)}\\
&+\frac{x_1 t^2 \left[p(x_3) x_2 (x_1-x_2) - p(x_2) x_3 (x_1-x_3) \right] }{(x_1-x_2) (x_1-x_3)(x_2-x_3) (1-x_3 t)}.
\end{split}
\ee


\subsection{Extension of the gate/wicket}
The operator which extends the gate or the wicket takes as input a formal power series with three catalytic variables and returns a formal power series of the same kind.

\begin{figure}[h!]
\psfrag{b1}{$b_1$}
\psfrag{b2}{$b_2$}
\psfrag{b3}{$b_3$}
\psfrag{c1}{$c_1$}
\psfrag{c2}{$c_2$}
\psfrag{c3}{$c_3$}
\includegraphics[width=6cm]{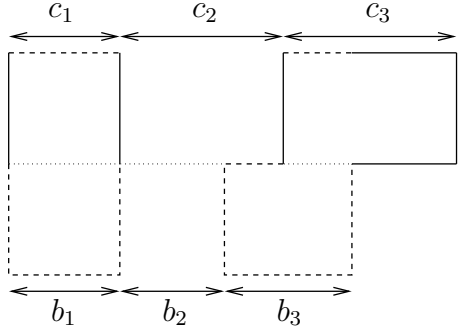}
\caption{The umbral evolution operator for extending the gate or wicket. The darkened lines contribute to the perimeter.} \label{fig:duringpunc}
\end{figure}

The nature of the \fd forces $c_1 = b_1$. The figure shows that $c_2 \geq b_2$ and $c_2 < b_2+b_3$. Thus the umbral operator $U_2$ acts as
\be
\begin{split}
x_1^{c_1} x_2^{b_2} x_3^{b_3} \mapsto &x_1^{c_1} \sum_{c_2 = b_2}^{b_2+b_3-1} \sum_{c_3 = b_2+b_3-c_2}^\infty x_2^{c_2} x_3^{c_3} t^{2+c_2+c_3-b_2-b_3} \\
& = \frac{t^2 x_1^{c_1} x_2^{b_2} x_3 (x_2^{b_3}-x_3^{b_3})}{(1-x_3 t)(x_2-x_3)}
\end{split}
\ee

Thus, $U_2$ acts on formal power series according to
\be
U_2(p(x_1,x_2,x_3)) = \frac{x_3 t^2}{(1-x_3 t)(x_2-x_3)} [ p(x_1,x_2,x_2) - p(x_1,x_2,x_3)]
\ee

Since Gated Ferrers diagrams are essentially unfinished Wicketed \fds we are in a position to describe the former. We can now express the generating function for these objects completely.

\begin{theo}
Let $\phi(t)$ be the half-perimeter generating function of \\Gated Ferrers diagrams. Then $\phi$ satisfies the following quadratic equation:
\be \label{punceq}
 \left[t^2 (1- 2 t)^4 (1-3t + t^2)\right] \phi^2 + \left[t^4 (1-3t+t^2) (1-2 t)^2 \right] \phi -t^{10}= 0 
\ee
\end{theo}

Before going on to the proof, it might be helpful to say a few words on how one could arrive at this answer. Using the umbral transfer matrix method and sufficient computing time, one simply generates polynomials in $x_1,x_2,x_3$ which correspond to a fixed half-perimeter value. Once one has a sufficient number of terms, one checks for a P-recursive ansatz.

In this case, there are three symbolic variables and therefore both generating of the data as well as checking the ansatz takes very long and therefore, one plugs in various integer values of the $x_i$'s and finds that for all of them, the equation satisfied by the conjectured generating function factorizes beautifully. Then a simple polynomial interpolation algorithm does the trick.

As far as we can tell, the sequence enumerated by $\phi$ have not been studied before. Being responsible mathizens, we entered it in {\em the} database as A133106 \cite{sloane}.

\begin{proof}
Let $\phi_{123}(x_1,x_2,x_3,t)$ be the generating function for Gated \fds where $t$ counts the half-perimeter, $x_1$ counts the width before the start of the wicket, $x_2$ counts the width of the wicket and $x_3$ counts the width after the end of the wicket. We claim that \\ $\phi_{123}(x_1,x_2,x_3,t)$ satisfies the equation below.

\be \label{punceq123}
\begin{split}
& \left[(x_2-x_3+x_3 (x_3-x_2) t+x_3 t^2) (1-x_1 t- t)^2 (1-x_2 t-2t + t^2)  \right.\\
&\left. \times (1-x_3 t- t)^2 \right] \phi_{123}^2 + \left[t^4 (x_1 x_2 x_3) (1-x_2 t -2t+t^2) (1-x_3 t+ t) \right. \\
& \left. \times (1-x_3 t - t) (1-x_1 t- t) \right] \phi_{123} -t^{10} (x_1 x_2 x_3)^2= 0 
\end{split}
\ee

To see that, notice that $\phi_{123}(x_1,x_2,x_3,t)$  satisfies the umbral evolution equation
\be \label{puncumb}
\phi_{123} = U_1 \left(\frac{x t^2}{1-t-xt} \right) + U_2(\phi_{123})
\ee
because $U_1$ takes a \fd and begins the gate which can be extended using $U_2$.

One solves \eqref{punceq123} for $\phi$ and notices that only the solution with the positive root has the correct Taylor expansion. Then one simply plugs it into \eqref{puncumb} and finds that it is satisfied. 

Now substitute $x_1=x_2=x_3=1$ into \eqref{punceq123} to find that $\phi_{123}(1,1,1,t)\\=\phi(x,t)$ satisfies \eqref{punceq}.
\end{proof}

\subsection{Termination of the wicket}
To describe Wicketed \fds, we also need the umbral operator describing the end of the wicket. This is essentially the reverse process of adding the wicket. Now we go from three catalytic variables to one.

\begin{figure}[h!]
\psfrag{a}{$a$}
\psfrag{b1}{$b_1$}
\psfrag{b2}{$b_2$}
\psfrag{b3}{$b_3$}
\includegraphics[width=6cm]{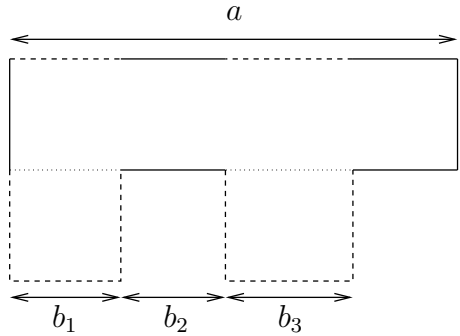}
\caption{The umbral evolution operator for ending the wicket. The darkened lines contribute to the perimeter.} \label{fig:endpunc}
\end{figure}

The action of $U_3$ is relatively straightforward.
\be
\begin{split}
x_1^{b_1} x_2^{b_2} x_3^{b_3} \mapsto & \sum_{a=b_1+b_2+b_3}^\infty x^a t^{1+a-b_1-b_3} \\
& = \frac{x^{b_1+b_2+b_3} t^{1+b_2}}{1-t x}
\end{split}
\ee

Thus, $U_3$ acts on formal power series by
\be \label{endpunc}
U_3(p(x_1,x_2,x_3)) = \frac{t}{1-t x} p(x,t x,x)
\ee

Since one needs to start with a Gated \fd to terminate the wicket, the initial condition for the evolution of \fds with one \fp will be $U_3(\phi)$.

\subsection{Termination of the diagram}
Once the wicket has been sealed up, the umbral operator is again given by \eqref{outpunc}. This is because the evolution of the \fd will continue as usual once the wicket has been established.

\begin{theo}
Let $\psi(t)$ be the half-perimeter generating function of Wicketed Ferrers diagrams. Then $\psi$ satisfies the following quadratic equation
\be \label{fdfpeq}
(2 t-1)^8 \psi^2 -t^6 (2 t-1)^4 \psi + t^{14} = 0.
\ee
\end{theo}

As far as we can tell, the sequence enumerated by $\psi$ has also not been studied before. We entered the sequence enumerated by $\psi$ as A133107  \cite{sloane}.

\begin{proof}
Let $\psi_1(x,t)$ be the generating function of Wicketed \fds where $t$ counts the half-perimeter and $x$ counts the width of the topmost segment. We claim that $\psi_1$ satisfies the equation below.

\be \label{fdfpeq1}
\begin{split}
&\left[(x t+ t-1)^6 (x t^2- t^2+2 t-1) (x t^2- t^2-x t- t+1) \right] \psi_1^2 \\
& - \left[t^6 x^2 (xt+ t-1)^3 (x t- t-1) (x t^2-t^2+2 t-1) \right] \psi_1 -t^{14} x^5= 0.
\end{split}
\ee

To see this, notice that $\psi_1$ satisfies the umbral equation \eqref{fdfpeq1} where the initial condition consists of \fds with a just-finished \fp and the evolution is the usual \fd evolution given by \eqref{outpunc}. That is,
\be \label{fdfpumb}
\psi_1 = U_3(\phi_{123}) + U_0(\psi_1)
\ee

One solves the quadratic equation \eqref{fdfpeq1} for $\psi_1$ and takes the negative root, which is the only one that has the correct Taylor expansion. One plugs the solutions of $\phi_{123}$ and $\psi_1$ into \eqref{fdfpumb} and uses the definition of the operators in \eqref{outpunc},\eqref{endpunc} to easily check that it is verified.

Finally, one substitutes $x=1$ into \eqref{fdfpeq1} to find that $\psi_1(1,t) = \psi(t)$ satisfies \eqref{fdfpeq}.
\end{proof}

\section{Remarks}
The generating function for Wicketed \fds is very pretty. Solving for $\psi$ in \eqref{fdfpeq} yields
\be
\psi(t) = t^6 \frac{1 - \sqrt{1- 4 t^2}}{2 (1-2 t)^4 }
\ee
and after some factorizing gives
\be
\psi(t)=\frac{1 - \sqrt{1- 4 t^2}}{2 t^2} \left( \frac{t^2}{1-2t} \right)^4.
\ee
which is a product of the half-perimeter generating function of staircase polygons (which gives the Catalan numbers) and \emph{four} sets of Ferrers diagrams as shown in \eqref{fdgf}. 

This nice factorization suggests an alternative proof using the technique of squeezing \cite{zeilcom}. We give some details of the idea in Appendix \ref{sec:altpf}.

One might wonder whether all the machinery used here is superfluous and whether this generating function might be arrived at by simpler means. We do not believe it to be so. We present a couple of arguments in favour of this assertion.

First off, a direct bijection is not so trivial because the first two terms in the series expansion of $\psi$ are 1 and 8 and in both of them the Catalan part is trivial (i.e. contributes 1). The series grows sufficiently fast so as to make experimentation difficult. Any natural bijection would be extremely interesting (look at Section \ref{sec:prize}).

Secondly, the Gessel-Viennot determinant formula \cite{gv} is not applicable here because the beginning and end of the Ferrers wicket are not constrained to lie at the corners of the \fd and in fact, are forbidden from being there. 

\section{Prize} \label{sec:prize}
We would like to offer a small award for a {\em nice, natural} bijection between the set of Wicketed Ferrers diagrams and the union of four sets of Ferrers diagrams and the set of staircase polygons\footnote{Any other set counted by the Catalan numbers would do just as well.}. The award consists of 50 of whatever basic unit of currency the winner decides (50 kuwaiti dinars or 50 rupees would be okay, 50 100-euro notes would not). The prize money will be paid in possibly an alternate currency to the winner using the exchange rate on the day the bijection is verified to be correct.

\section{Acknowledgements}
The author would like to express his heartfelt thanks to Doron Zeilberger for patiently explaining the umbral transfer matrix method, for many discussions, and for reading and correcting earlier versions of the draft. The author would also like to thank Ira Gessel for a discussion about Gessel-Viennot.

\appendix

\section{Sketch of an Alternate Proof} \label{sec:altpf}
The idea is to consider the generating function of a pair of nonintersecting paths satisfying some constaints imposed by the structure of the Wicketed Ferrers Diagrams. These are enumerated by the sum of their lengths, which is essentially the half-perimeter. 

Note that one can add horizontal and vertical segments, both just after the first segment as well as just before the last segment to get a new legal pair of paths. This can be done both for the inner and the outer path. Doing it for the outer path properly gives a factor of $(t/(1-2t))^2$. 

In an ideal world the same could be done for the inner path to give the same factor again. What would then be left is precisely the usual staircase polygons, which give the Catalan generating function. The problem is that various constraints come into play for the inner path because the distance between both of them have to be atleast one both horizontally and vertically. 

\begin{figure}[h!]
\psfrag{a}{$\alpha$}
\psfrag{b}{$\beta$}
\includegraphics[width=6cm]{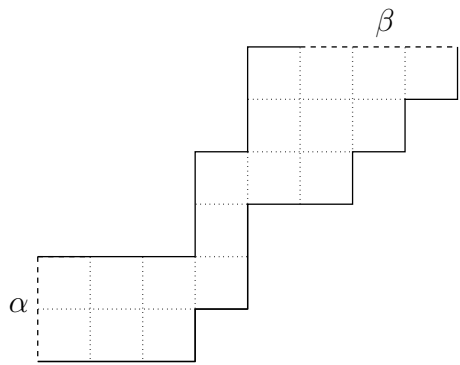}
\caption{A nibbled staircase polygon contributing $\alpha^2 \beta^3$. Only the darkened lines contribute to the perimeter.} \label{fig:nibsp}
\end{figure}

At this point one has reduced the problem to that of calculating the half-perimeter generating function of {\em Nibbled Staircase Polygons}  with nibbling parameters $\alpha$ and $\beta$ like that in Figure \ref{fig:nibsp}, where the walk on top is the inner path and the one on the bottom is the outer path. The idea is that some part of the first vertical segment and the last horizontal segment of the higher path does not contribute to the half-perimeter. The generating function of these nibbled staircase polygons can be calculated again using the umbral transfer matrix method.

\end{document}